\documentclass[12pt]{amsart}
\usepackage{amsfonts,amssymb,amscd,amsmath,amsthm,enumerate,verbatim,calc}

\theoremstyle{plain}
\newtheorem{Theorem}{Theorem}[section]

\newtheorem{Lemma}[Theorem]{Lemma}
\newtheorem{Corollary}[Theorem]{Corollary} 
\newtheorem{Proposition}[Theorem]{Proposition}
\newtheorem*{Question}{Question}

\newtheorem*{Conjecture}{Conjecture}
\theoremstyle{definition}
\newtheorem{Definition}[Theorem]{Definition}
\newtheorem{Example}[Theorem]{Example}
\newtheorem{Remark}[Theorem]{Remark}

\theoremstyle{remark}

\newtheorem*{chunk*}{}

\numberwithin{equation}{Theorem}

\newcommand{\mlabel}[1]%
  {\mbox{}\marginpar{\raggedleft\hspace{0pt}{\rm\ttfamily#1}}\label{#1}}

\newcommand{\fm}{{\mathfrak m}}

\newcounter{hours}\newcounter{minutes}

\newcommand{\excise}[1]{}
\begin{document}
\title[The Frobenius map on local cohomology modules in mixed characteristic]
{The Frobenius action on local cohomology modules in mixed characteristic}
\author[K.Shimomoto]{Kazuma Shimomoto}
\address{School of Mathematics, University of Minnesota, 127 Vincent Hall, 206 Church St. S.E.,
Minneapolis, MN 55455 USA.}
\email{shimo027@umn.edu}


\begin{abstract}
R. Heitmann's proof of the Direct Summand Conjecture has opened a new approach to the study of homological conjectures in mixed characteristic. Inspired by his work and by the methods of almost ring theory, we discuss a normalized length for certain torsion modules, which was introduced by G. Faltings. Using the normalized length and the Frobenius map, we prove some results of local cohomology for local rings in mixed characteristic, which has an immediate implication for the subject of splinters studied by A. Singh. 
\end{abstract}

\maketitle

\vspace*{6pt}

\section{Introduction}
\label{sec:introduction}

In the present paper, we apply the Frobenius map to prove some results on local cohomology modules of local rings of mixed characteristic. These results were motivated by connections to the homological conjectures, in particular the Direct Summand Conjecture which states the following:

\begin{Conjecture}
Let $R$ be a regular local ring and let $R \to S$ be a module-finite extension. Then $R$ is a direct summand of $S$ as an $R$-module.
\end{Conjecture}

Let $(R,\fm)$ be a Noetherian local ring. Recall that an $R$-algebra $T$ satisfying $\fm T \ne T$ is a ($balanced$) $big$ $Cohen$-$Macaulay$ $R$-$algebra$ if every system of parameters of $R$ is a regular sequence on $T$. If $T$ is a big Cohen-Macaulay algebra, then the local cohomology modules $H^{i}_{\fm}(T)$ are zero for all $i < \mathrm{dim}~R$. The importance of the notion of such algebras is contained in the following (see~\cite{H2}):

\begin{Proposition}
Let $(R,\fm)$ be a complete local domain of arbitrary characteristic. If there exists a big Cohen-Macaulay $R$-algebra, then the Direct Summand Conjecture holds.
\end{Proposition}

Let $R^{+}$ be the integral closure of $R$ in the algebraic closure of the field of fractions of $R$. Then it is shown in \cite{HH1} that $R^{+}$ is a big Cohen-Macaulay $R$-algebra if $R$ has characteristic $p>0$. In the mixed characteristic case, the following result was recently established by Heitmann (see~\cite{He},~\cite{He2}) to prove the Direct Summand Conjecture in dimension 3:

\begin{Theorem}
Let $(R,\fm)$ be a 3-dimensional complete local domain of mixed characteristic $p>0$. Then $c^{\epsilon} \cdot H^{2}_{\fm}(R^{+})=0$ for any $c \in \fm$ and $\epsilon \in \mathbb{Q}$, $\epsilon>0$. 
\end{Theorem}

As a corollary, Hochster (see~\cite{H3}) deduced the existence of big Cohen-Macaulay algebras:

\begin{Theorem}
Let $(R,\fm)$ be the same as in the above theorem. Then there exists a big Cohen-Macaulay $R$-algebra in a weakly functorial sense. 
\end{Theorem}

In connection with the above result, an interesting question is whether $R^{+}$ is ``almost'' Cohen-Macaulay in higher dimension in the sense that the local cohomology modules $H^{i}_{\fm}(R^{+})$ are ``$v$-almost zero'' (we give a definition of this term using a valuation in Definition 2.7) for all $i < \mathrm{dim}~R$. In fact, once this is answered positively, Heitmann's idea works everywhere to prove the Direct Summand Conjecture in any dimension. We also remark that the Direct Summand Conjecture is a theorem for all equicharacteristic rings (see~\cite{H1} for the proof), while it is an open question in dimension $>3$ in mixed characteristic. Now let $(R,\fm)$ be a $d$-dimensional power series ring over a complete discrete valuation ring $V$ of mixed characteristic $p>0$ with perfect residue field. Let $R \to S$ be a module-finite extension of local domains. By Cohen structure theorem, every complete local ring is always module-finite over some power series ring. In particular, $R^{+}=S^{+}$.

The structure of this paper goes as follows. First we define a certain big ring $R_{\infty}$ that is obtained as a flat colimit of $R$ and then introduce a notion of $normalized$ $length$ $\lambda_{\infty}(M) \in \mathbb{R} \cup \{ \infty \}$
for an $\fm$-torsion $R_{\infty}$-module $M$, which was introduced by Faltings in his method of almost $\mathrm{\acute{e}}$tale extensions (see~\cite{F}). Then the aim of this paper is to investigate the following question: 

\begin{Question}[1]
Let $\Phi:H^{i}_{\fm}(S)\otimes_{R} R_{\infty} \to H^{i}_{\fm}(S^{+})$ be the map of local cohomology modules induced by the natural map $S \otimes_{R} R_{\infty} \to S^{+}$. Then is it true that $\lambda_{\infty}(\mathrm{Im}(\Phi))=0$ for all $i<d$? 
\end{Question}

Faltings investigated this question when the localization map $R[p^{-1}] \to S[p^{-1}]$ is $\mathrm{\acute{e}}$tale. Our main result (see Theorem 3.6) gives a partial answer to this question without any assumptions on the module-finite extension $R \to S$. Question 1 is related to a notion of almost zero modules (see~\cite{GR1} for the definition). We shall define a notion of $v$-almost zero modules below, which suffices for the study of local cohomology modules. The importance of Question 1 is contained in the following:

\begin{Proposition}
If Question 1 has a positive answer for every module-finite extension $R \to S$, then $R^{+}$ is almost Cohen-Macaulay. In particular, the Direct Summand Conjecture holds.
\end{Proposition}

The proof of this proposition follows from Proposition 2.15 together with Heitmann's idea to deduce the Direct Summand Conjecture from the almost vanishing of local cohomology modules. Proposition 2.15 is essential in studying $v$-almost zero modules. The main point in working with $S^{+}$ is that the Frobenius endomorphism on $S^{+}/pS^{+}$ is surjective. Hence it induces a unique ring isomorphism $\overline{F}_{S^{+}}:S^{+}/p^{1/p}S^{+} \simeq S^{+}/pS^{+}$ and 
$$
H^{i}_{\fm}(S^{+}/p^{1/p}S^{+}) \simeq H^{i}_{\fm}(S^{+}/pS^{+})^{[F]},
$$
in which case the right hand side is viewed as an $R_{\infty}$-module via the Frobenius map. Then one applies Theorem 2.12 to study the above local cohomology map. To obtain non-trivial results, we need to restrict our attention to the finite length cohomology $H^{k}_{\fm}(S)$ together with an additional assumption that $H^{k-1}_{\fm}(S)$ is zero, in which case $\Phi:H^{k}_{\fm}(S) \otimes_{R} R_{\infty} \to H^{k}_{\fm}(S^{+})$ is, at least, shown to be not injective if $H^{k}_{\fm}(S)$ is non-zero.

The second main theorem in this paper is to state a certain type of finiteness conditions under which the local cohomology modules $H^{i}_{\fm}(S_{\infty})$ are $v$-almost zero for $i<d$ for some big normal domain $S_{\infty}$ of mixed characteristic that contains both $S$ and $R_{\infty}$ and that satisfies the condition that the Frobenius endomorphism
$$
F_{S_{\infty}}:S_{\infty}/pS_{\infty} \to S_{\infty}/pS_{\infty} 
$$
is surjective. We use normalized length to produce such conditions (see Theorem 3.4). At this point, as it is not clear as to the class of mixed characteristic rings on which the Frobenius map is surjective after killing $p$, we leave this issue in this paper. Finally, we point out that the idea of Faltings' work stems from an attempt to extend the classical Nagata-Zariski's purity theorem to its almost analogue.

\section{Normalized length}
\label{sec:normalized-length}

Throughout this article, $(R,\fm,k)$ will be a complete regular local ring of mixed characteristic $p>0$ $$ 
V[[x_2,\ldots,x_d]] 
$$ 
such that $(V,pV,k)$ is a complete discrete valuation ring with perfect residue field $k$ of characteristic $p>0$. Then we consider the direct system of regular local rings
$$
R=R_0 \hookrightarrow R_1 \hookrightarrow \cdots \hookrightarrow R_e \hookrightarrow \cdots
$$ 
where $R_e := V[p^{1/p^e}][[x_2^{1/p^e},\dots,x_d^{1/p^e}]]$ with maximal ideal $\fm_{e}$, and we let $R_{\infty}:=\bigcup_{e \in \mathbb{N}} R_{e}$. Then the ring $R_{\infty}$ possesses the following properties:

\begin{enumerate}
\item
$R_{\infty}$ is a henselian quasilocal ring of dimension $d$ with unique maximal ideal $\fm_{\infty}$. 

\item
$R_{\infty}$ is faithfully flat and integral over $R$.

\item
The Frobenius endomorphism $F_{R_{\infty}}:R_{\infty}/pR_{\infty} \to R_{\infty}/pR_{\infty}$
is surjective.
\end{enumerate}

To see this, note that the first two properties are stable under infinite ascending unions of rings, while the third follows from the surjectivity: $R_{e+1}/pR_{e+1} \twoheadrightarrow R_{e}/pR_{e}$ induced by the Frobenius map. For more on henselian local rings, see~\cite{N}.

Let $\mathcal{M}$ denote the category of $\fm$-torsion $R_{\infty}$-modules where $\fm$ is the maximal ideal of $R$. An $R_{\infty}$-module $M$ is said to be $\fm$-$torsion$ if every element of $M$ is annihilated by some power of $\fm$. Examples are local cohomology modules. Let us start with the definition of normalized length. 

\begin{Definition} \

\begin{enumerate}

\item
Let $M \in \mathcal{M}$ be a finitely presented module. Then there is a finite presentation for $M$:
$$
\begin{CD}
R^{\oplus m_{1}}_{\infty} @>\varphi>> R^{\oplus m_{0}}_{\infty} @>>> M @>>> 0 \cdots(*). 
\end{CD}
$$
For sufficiently large $n \in \mathbb{N}$, we have $a_{ij} \in R_{n}$ with $\varphi=(a_{ij})$, and for some finite $R_{n}$-module $M_{n}$ we have 
$$
\begin{CD}
R^{\oplus m_{1}}_{n} @>\varphi>> R^{\oplus m_{0}}_{n} @>>> M_{n} @>>> 0 \cdots(**). \\
\end{CD}
$$
Since $R_{\infty}$ is flat over $R_{n}$, we may tensor $R_{\infty}$ with $(**)$ over $R_{n}$ to obtain
$$
\begin{CD}
R^{\oplus m_{1}}_{\infty} @>\varphi>> R^{\oplus m_{0}}_{\infty} @>>> R_{\infty} \otimes_{R_n} M_{n} @>>> 0. \\
\end{CD}
$$
which is equivalent to the presentation $(*)$, and thus $M \simeq R_{\infty} \otimes_{R_n} M_{n}$. Then define
$$
\lambda_{\infty}(M):=\frac{1}{p^{dn}} \cdot \ell(M_{n}) \in \mathbb{R}.
$$

\item
Let $M \in \mathcal{M}$ be a finitely generated module. Then define
$$
\lambda_{\infty}{(M)}:=\inf_{N \twoheadrightarrow M} \lambda_{\infty}(N) \in \mathbb{R}
$$
where ``inf'' is taken over all finitely presented modules $N$ in $\mathcal{M}$ that map onto $M$.

\item
Let $M \in \mathcal{M}$ be an arbitrary module. Then define
$$
\lambda_{\infty}(M):=\sup_{N \subset M} \lambda_{\infty}(N) \in \mathbb{R} \cup \{\infty\}
$$
where ``sup'' is taken over all finitely generated submodules $N$ of $M$.
\end{enumerate}

\end{Definition}

\begin{Lemma}
The normalized length $\lambda_{\infty}(M)$ is well-defined for $M \in \mathcal{M}$. 
\end{Lemma}

\begin{proof}
Let $M \in \mathcal{M}$ be a finitely presented module. Then one has to verify that $\lambda_{\infty}(M)$ is independent of $n$ in the definition. Let $m \ge n$. Then since $R_{n} \to R_{m}$ is flat and $M_{m} \simeq R_{m} \otimes_{R_n} M_{n}$, we have:
$$
\lambda_{\infty}(M)=\frac{1}{p^{dn}} \cdot \ell(M_{n})=\frac{1}{p^{dm}} \cdot \frac{1}{p^{d(n-m)}} \cdot \ell(M_{n})=\frac{1}{p^{dm}} \cdot \ell(R_{m} \otimes_{R_n} M_{n})=\frac{1}{p^{dm}} \cdot \ell(M_{m}),
$$
which is the claim.

Next, one has to verify that ii) coincides with i) for finitely presented $M \in \mathcal{M}$. Let $N \in \mathcal{M}$ be a finitely presented module that maps onto $M$. Then it suffices to see that $\lambda_{\infty}(N) \ge \lambda_{\infty}(M)$ since $M$ maps onto itself. This will be seen by taking a finite presentation $N_{n} \twoheadrightarrow M_{n}$ for $n \gg 0$ and hence $\ell(N_{n}) \ge \ell(M_{n})$.

Finally, one needs to verify that iii) coincides with ii) for finitely generated $M \in \mathcal{M}$. Let $N \subset M$ be any finitely generated submodule. Let us choose finite sets of generators $\Sigma \subset N$ and $\Sigma' \subset M$
such that $\Sigma \subset \Sigma'$. If $N_{n}$ and $M_{n}$ denote the $R_{n}$-submodules generated by $\Sigma$ and $\Sigma'$, respectively, then $\ell(N_{n}) \le \ell(M_{n})$ for all $n \in \mathbb{N}$. Hence the desired claim follows from Lemma 2.3 below.
\end{proof}

We used the following lemma in the proof of Lemma 2.2.

\begin{Lemma}
Let $M \in \mathcal{M}$ be a finitely generated module and let $\Sigma$ be any system of generators of $M$. Let $M_{n}$ denote the $R_{n}$-module generated by $\Sigma$. Then $M_{n} \otimes_{R_n} R_{\infty} \in \mathcal{M}$ is a finitely presented module that admits a surjection onto $M$, and
$$
 \lambda_{\infty}(M)=\lim_{n \to \infty}\lambda_{\infty}(M_{n} \otimes_{R_n} R_{\infty})=\lim_{n \to \infty} \frac{1}{p^{dn}}\cdot \ell(M_{n}).
$$
\end{Lemma}

\begin{proof}
Let $\Sigma$ and $\Sigma'$ be finite sets of generators of the module $M$. We denote by $M_{n}$ and $M_{n}'$ the $R_{n}$-modules generated by $\Sigma$ and $\Sigma'$, respectively. Then it is easy to see that $M_{n}=M_{n}'$ for sufficiently large $n>0$. So $\ell(M_{n})=\ell(M_{n}')$ for $n \gg 0$. By definition, it suffices to prove if $N$ is a finitely presented module in $\mathcal{M}$ that surjects onto $M$, then there exists some $n \in \mathbb{N}$ such that the $R_{\infty}$-module $N$ admits a surjection onto $M_{n} \otimes_{R_n} R_{\infty}$ as we would get 
$$
\lambda_{\infty}(M) \le \lambda_{\infty}(M_{n} \otimes_{R_n} R_{\infty}) \le \lambda_{\infty}(N),
$$
from which the conclusion easily follows.

Since $N$ is finitely presented, there exists some $n \in \mathbb{N}$ such that $N \simeq N_{n} \otimes_{R_n} R_{\infty}$ for some finite $R_{n}$-module $N_{n}$. Then define $M_{n}$ to be image of the composite map $N_{n} \hookrightarrow N_{n} \otimes_{R_n} R_{\infty} \twoheadrightarrow M$ and the module $N$ obviously admits a surjection onto $M_{n} \otimes_{R_n} R_{\infty}$.
\end{proof}

The following lemma follows from the part iii) of the definition.

\begin{Lemma}
Let $M \in \mathcal{M}$ and let $N \subset M$ be its submodule. Then $\lambda_{\infty}(N) \le \lambda_{\infty}(M)$.
\end{Lemma}

\begin{Remark}
One can extend the class of rings for which the normalized length is well-defined as follows. For the definition of the length, we only use the fact that $R_{e} \hookrightarrow R_{e+1}$ is flat and of the constant rank $p^{d}$. It is therefore natural to start with a family of flat extensions of Noetherian local rings
$$
R:=R_{0} \subset\ \cdots \subset R_{n} \subset R_{n+1} \subset \cdots \subset R_{\infty}
$$ 
such that the rank $[R_{n+1}:R_{n}]$ stabilizes for $n \gg 0$ and $\sqrt{\fm_{n}R_{n+1}}=\fm_{n+1}$ for all $n \in \mathbb{N}$. Under this set-up, one can define the normalized length as in the case $R$ is a complete regular local ring. We refer the reader to \cite{GR2} for more detail, where the normalized length is discussed for a more general class of rings. 
\end{Remark}

\begin{Remark}
From the definition, $\lambda_{\infty}(M)$ is finite for any finitely generated $M \in \mathcal{M}$. It is also true that if $M$ is finitely presented in $\mathcal{M}$, then: 
$$
\lambda_{\infty}(M)=0 \iff M=0,
$$
which does not hold for finitely generated modules. For example, take $R_{\infty}/\fm_{\infty} \in \mathcal{M}$, which is not a finitely presented $R_{\infty}$-module, and it is easy to see that $\lambda_{\infty}(R_{\infty}/\fm_{\infty})=0$. 
\end{Remark}

This remark suggests the following definition.

\begin{Definition}
Let $M$ be an $R_{\infty}$-module in $\mathcal{M}$, and let $v$ be a valuation on $R_{\infty}$ that is positive on the maximal ideal of $R_{\infty}$. Then we say that
\begin{enumerate}
\item
$M$ has $almost$ $finite$ $length$ if for any rational number $\epsilon>0$, there exists an element $b \in R_{\infty}$ such that $\lambda_{\infty}(b \cdot M)< \infty$ and $v(b) < \epsilon$. 

\item 
$M$ is $v$-$almost$ $zero$ if for any $m \in M$ and rational number $\epsilon>0$, there exists an element $b \in R_{\infty}$ such that $b \cdot m=0$ and $v(b) < \epsilon$. 
\end{enumerate}
\end{Definition}

In Proposition 2.15, which is proved only for an $\fm$-adic valuation, we will relate the $v$-almost zero modules to the vanishing of the normalized length. The normalized length behaves additively on the short exact sequence of $R_{\infty}$-modules in $\mathcal{M}$.

\begin{Proposition}
Let $0 \to M' \to M \to M'' \to 0$ be an exact sequence of $\fm$-torsion $R_{\infty}$-modules. Then
$$
\lambda_{\infty}(M)=\lambda_{\infty}(M')+\lambda_{\infty}(M'').
$$
\end{Proposition}

\begin{proof}
In the case where all relevant modules are either finitely presented or finitely generated, the claim follows easily by the definition of the normalized length together with Lemma 2.3. So let us consider the general case.  Since the image of every finitely generated submodule of $M$ in $M''$ is finitely generated, it suffices to consider the case in which both $M$ and $M''$ are finitely generated by the definition of normalized length for arbitrary modules. Let us denote by $M_{n}$ an $R_{n}$-submodule of $M$ generated by any fixed generators of the $R_{\infty}$-module $M$. Let $M'_{n}:=M' \cap M_{n}$ and let $M''_{n}$ be the image of $M_{n}$ in $M''$. Then we have a short exact sequence: $0 \to M_{n}' \to  M_{n} \to M_{n}'' \to 0$.

For a finite set of elements of $M'$ the $R_n$-submodule $\widetilde{M}_{n}'$ generated by them is contained in some $M_{n}'$ and satisfies that $\ell(\widetilde{M}_{n}') \le \ell(M_n')$ for $n \gg 0$. Then since $M$ and $M''$ are finitely generated, it follows from Lemma 2.3 that 
$$
\lambda_{\infty}(M')=\lim_{n \to \infty} \limsup_{\widetilde{M}_{n}' \subset M'} \frac{1}{p^{dn}}\ell(\widetilde{M}_{n}') \le \lim_{n \to \infty}\frac{1}{p^{dn}}\ell(M_{n}')=\lambda_{\infty}(M)-\lambda_{\infty}(M''),
$$
which is equivalent to the inequality
$$
\lambda_{\infty}(M)-\lambda_{\infty}(M'') \ge \lambda_{\infty}(M').
$$
On the other hand, since $M''=\varinjlim_{n} M_{n}''\otimes_{R_n} R_{\infty}$, for any $\epsilon>0$ there exists sufficiently large $n>0$ such that
$$
\lambda_{\infty}(M'') \ge \lambda_{\infty}(M_{n}''\otimes_{R_n} R_{\infty})-\epsilon.
$$
Let us consider the following commutative diagram with exact rows.
$$
\begin{CD}
0 @>>> M_{n}' \otimes_{R_{n}} R_{n+m} @>>> M_{n} \otimes_{R_{n}} R_{n+m} @>>> M_{n}'' \otimes_{R_{n}} R_{n+m} @>>> 0 \\
@. @V f VV @V g VV @V h VV \\
0 @>>> M_{n+m}' @>>> M_{n+m} @>>> M_{n+m}'' @>>> 0 \\
\end{CD}
$$
Since $\mathrm{Coker}(g)=0$, there is an epimorphism $\mathrm{Ker}(h) \twoheadrightarrow \mathrm{Coker}(f)=M_{n+m}'/R_{n+m}\cdot M_{n}'$, which follows from the snake lemma. Let $\gamma=\ell(\mathrm{Ker}(h))$. Then there follows that
$$
\gamma=\ell(M_{n}''\otimes_{R_n} R_{n+m})-\ell(M_{n+m}'')=p^{dm} \cdot \ell(M_{n}'')-\ell(M_{n+m}'').
$$
It also follows by our assumption that
$$
\lambda_{\infty}(M'') \ge \frac{\ell(M_{n}'')}{p^{dn}}-\epsilon.
$$
We can put all inequalities obtained above together to obtain that
$$
\frac{\gamma}{p^{d{(n+m)}}} \le \lambda_{\infty}(M'')-\frac{\ell(M_{n+m}'')}{p^{d{(n+m)}}}+\epsilon \le \epsilon.
$$
On the other hand, since $\{\lambda_{\infty}(M_{n+m}''\otimes_{R_{n+m}} R_{\infty})\}_{m\in\mathbb{N}}$ is a decreasing sequence in $\mathbb{Q}$, we get the inequality $\lambda_{\infty}(M'') \le \lambda_{\infty}(M_{n+m}''\otimes_{R_{n+m}} R_{\infty})$ for all $m \in \mathbb{N}$. Now we have $\gamma \le \epsilon \cdot p^{d{(n+m)}}$.  This inequality together with an epimorphism $\mathrm{Ker}(h) \twoheadrightarrow \mathrm{Coker}(f)$ yields that
$$
\ell(M_{n+m}') \le \ell(R_{n+m}\cdot M_{n}')+\epsilon \cdot p^{d{(n+m)}}
$$
and therefore
$$
\lambda_{\infty}(M)-\lambda_{\infty}(M'') \le \lim_{n \to \infty}\frac{1}{p^{dn}} \cdot \ell(M_{n}') \le \lim_{m \to \infty}\frac{1}{p^{d{(n+m)}}} \cdot \ell(R_{n+m} \cdot M_{n}')+\epsilon
$$
$$
=\lambda_{\infty}(R_{\infty}\cdot M_{n}')+\epsilon \le \lambda_{\infty}(M')+\epsilon.
$$
Since $\epsilon$ can be taken arbitrarily small, $\lambda_{\infty}(M)-\lambda_{\infty}(M'') \le \lambda_{\infty}(M')$, and hence the claim follows.
\end{proof}

\begin{Proposition}
Let $N$ be an $R_{\infty}$-module that admits an injection into some finitely presented module $M \in \mathcal{M}$. If $\lambda_{\infty}(N)=0$, then $N=0$. 
\end{Proposition}

\begin{proof}
As $N=0$ if and only if every finitely generated submodule of $N$ is zero, we may harmlessly assume $N$ is finitely generated. Let us consider the exact sequence
$$
0 \to N \to M \to M/N \to 0.
$$
Then the quotient module $M/N$ is finitely presented with $\lambda_{\infty}(M)=\lambda_{\infty}(M/N)$. However, since $M$ is finitely presented, there exists an $R_{n}$-submodule $M_{n} \subset M$ such that $R_{\infty} \otimes_{R_n} M_{n} \simeq M$. Let $M_{n}'$ be the image of $M_{n}$ under the surjection $M \twoheadrightarrow M/N$. Then it can be shown that $R_{\infty} \otimes_{R_n} M_{n}' \simeq M/N$ by computing the finite presentation of $M/N$. Then we must have $\ell(M_{n})=\ell(M_{n}')$, which is the case only when $M_{n}=M_{n}'$ and $M=M/N$. Hence $N=0$, which is the desired result.
\end{proof}

\begin{Definition}
Let $M$ be any $R_{\infty}$-module such that $p \cdot M=0$. Then $M^{[F]}$ denotes the $R_{\infty}$-module which is isomorphic to $M$ as an abelian group and has the module structure given by $a \cdot m:=a^{p}m=:m \cdot a$ for $a \in R_{\infty}$ and $m \in M^{[F]}$.  
\end{Definition}

Note that $M \mapsto M^{[F]}$ defines the identity functor on the underlying abelian groups and we give the bimodule structure on the modules via the Frobenius map.

\begin{Example}
One expects that it would be very useful if the like Frobenius maps could be defined for rings of mixed characteristic. Here is an example: Let $M$ be an $R_{\infty}/p^{1/p}R_{\infty}$-module that admits a Frobenius map $F$ in the following sense: the map $F:M \to M$ is a homomorphism of abelian groups. As an $R_{\infty}$-module map, let $F(am)=a^{p}F(m)$ for $a \in R_{\infty}$, $m \in M$, and hence $F((a+b)m)=(a+b)^{p}F(m)=a^{p}F(m)+b^{p}F(m)$. Now one has the following map: 
$$
M \otimes_{R_{\infty}} R_{\infty}^{[F]} \to M^{[F]}~:~m \otimes a \mapsto F(m)a,
$$
which is well-defined as $p^{1/p}$ annihilates the module $M$ by assumption and is called the $relative$ $Frobenius$ $map$. 
\end{Example}

We are now ready to prove the following ``Frobenius pull-back formula'', whose importance is expressed by studying certain torsion modules, in particular local cohomology modules, that admit a Frobenius action.

\begin{Theorem}
Let $M$ be an $\fm$-torsion $R_{\infty}/pR_{\infty}$-module. Then
$$   
\lambda_{\infty}(M^{[F]})=\frac{1}{p^d} \cdot \lambda_{\infty}(M).
$$
\end{Theorem}

\begin{proof}
Let us start with the case where $M$ is finitely presented in $\mathcal{M}$, and let
$$
\begin{CD}    
R_{\infty}^{\oplus m_{1}} @ >\varphi >> R_{\infty}^{\oplus m_{0}} @>>> M @ >>> 0 \\\end{CD}
$$
be the presentation for $M$ such that $\varphi$ is defined over $R_{n}$ for some $n \ge 0$. Since the functor $M \mapsto M^{[F]}$ is exact and $p \cdot M=0$ by assumption, there follows the commutative diagram
$$
\begin{CD}
(R_{\infty}/p^{1/p}R_{\infty})^{\oplus m_{0}} @ >>> (R_{\infty}/p^{1/p}R_{\infty})^{\oplus m_{1}}\\
@ V F^{m_{0}} V\wr V @ V F^{m_{1}} V\wr V \\
((R_{\infty}/pR_{\infty})^{[F]})^{ \oplus m_{0}} @ >>> ((R_{\infty}/pR_{\infty})^{[F]})^{\oplus m_{1}} @>>> M^{[F]} @ >>> 0 \\
\end{CD}
$$
where the vertical maps are induced by the ring isomorphism
$$
\overline{F}_{R_{\infty}}:R_{\infty}/p^{1/p}R_{\infty} \simeq (R_{\infty}/pR_{\infty})^{[F]}
$$
which is induced by the Frobenius map. Then using the above diagram together with the ring isomorphism $R_{n+1}/p^{1/p}R_{n+1} \simeq (R_{n}/pR_{n})^{[F]}$ induced by the Frobenius, we have the following exact sequence 
$$
\begin{CD}
(R_{n+1}/p^{1/p}R_{n+1})^{m_0} @>>> (R_{n+1}/p^{1/p}R_{n+1})^{m_1} @>>> M_{n}^{[F]} @>>> 0 \\
\end{CD}
$$
which in turn gives an $R_{\infty}$-module isomorphism $R_{\infty} \otimes_{R_{n+1}} M_{n}^{[F]} \simeq M^{[F]}$. On the other hand, since the residue class field of $R$ is assumed to be perfect, we have
$$
\ell(M_{n}^{[F]})=\ell(M_{n})
$$
and therefore 
$$
\lambda_{\infty}(M^{[F]})=\frac{1}{p^{d(n+1)}} \cdot \ell(M_{n}^{[F]})=\frac{1}{p^d} \cdot \frac{1}{p^{dn}} \cdot \ell(M_n)=\frac{\lambda_{\infty}(M)}{p^d},
$$
which completes the case of finitely presented modules.

Next let us assume $M$ is finitely generated in $\mathcal{M}$. Let $ N $ be a finitely presented module in $\mathcal{M}$ that admits a surjection $N \twoheadrightarrow M$.  Then this map factors as $N \twoheadrightarrow N/pN \twoheadrightarrow M$ by assumption and we deduce that
$$
\lambda_{\infty}(N) \ge \lambda_{\infty}(N/pN) \ge \lambda_{\infty}(M).
$$
Therefore we can replace $N$ with $N/pN$ to prove the theorem. Then since $N^{[F]}$ belongs to $\mathcal{M}$ and is a finitely presented module that maps onto $M^{[F]}$, we get $\lambda_{\infty}(N)=p^d\cdot\lambda_{\infty}(N^{[F]}) \ge p^{d}\cdot\lambda_{\infty}(M^{[F]})$ and thus $\lambda_{\infty}(M)\ge p^d\cdot\lambda_{\infty}(M^{[F]})$.

Conversely, let $N'$ be a finitely presented module that maps onto $M^{[F]}$. Then there is an $R_n$-module $N'_n$ such that $N'=R_{\infty} \otimes_{R_n} N'_n $. Now take $u_1,\ldots,u_s$ to be a set of generators of $M^{[F]}$ and define the $R_{n-1}$-module
$$
M_{n}:=R_{n-1}u_{1}+\cdots+R_{n-1}u_{s}\subset M^{[F]}.
$$
Then $M_n^{[F]} $ can be viewed as an $R_n$-module via the Frobenius $F:R_{n}/pR_{n} \to R_{n-1}/pR_{n-1}$, and the module $M_n^{[F]} $ is finitely generated over $R_{n} $. Hence we get the surjective map of $R_{\infty}$-modules as follows
$$
N'=R_{\infty} \otimes_{R_{n}} N'_{n} \twoheadrightarrow R_{\infty} \otimes_{R_{n}} M_n^{[F]} \twoheadrightarrow M^{[F]}.
$$
Now we claim that there is an $R_{\infty}$-module isomorphism:
$$
R_{\infty} \otimes_{R_{n}} M_{n}^{[F]} \simeq (R_{\infty} \otimes_{R_n} M_{n+1})^{[F]}.
$$
For a proof, let
$$
(R_{n-1}/pR_{n-1})^{\oplus m_0} \to (R_{n-1}/pR_{n-1})^{\oplus m_1} \to M_{n} \to 0
$$
be the presentation of the $R_{n-1}$-module $M_{n}$. Then the presentation of $M_{n+1}$ can be obtained by replacing $R_{n-1}$ by $R_n$ in the above presentation. Consider the following commutative diagram
{\small
$$
\begin{CD}
R_{\infty} \otimes_{R_n} ((R_{n-1}/pR_{n-1})^{[F]})^{\oplus m_0} @>>> R_{\infty} \otimes_{R_n} ((R_{n-1}/pR_{n-1})^{[F]})^{\oplus m_1} @>>> R_{\infty} \otimes_{R_n} M_{n}^{[F]} @ >>> 0 \\
@ V \Phi^{m_0} VV @ V \Phi^{m_1} VV @ VVV \\
(R_{\infty} \otimes_{R_{n}} (R_{n}/pR_{n})^{\oplus m_0})^{[F]} @ >>> (R_{\infty} \otimes_{R_{n}} (R_{n}/pR_{n})^{\oplus m_1})^{[F]} @ >>> (R_{\infty} \otimes_{R_{n}} M_{n+1})^{[F]} @ >>> 0 \\
@VV \wr V @VV \wr V  \\
((R_{\infty}/pR_{\infty})^{\oplus m_0})^{[F]} @ >>> ((R_{\infty}/pR_{\infty})^{\oplus m_1})^{[F]}\\
\end{CD}
$$
} 
in which both $\Phi^{m_0}$ and $\Phi^{m_1}$ are induced by the $ R_{\infty} $-module map 
$$
\Phi:R_{\infty} \otimes_{R_{n}}(R_{n-1}/pR_{n-1})^{[F]} \to (R_{\infty}/pR_{\infty})^{[F]}:~r \otimes m \mapsto r^{p}m
$$
and $(R_{n-1}/pR_{n-1})^{[F]}$ is viewed as an $R_n$-module via the Frobenius $F:R_{n}/pR_{n} \to R_{n-1}/pR_{n-1}$. Then $\Phi$ factors as
$$
R_{\infty} \otimes_{R_n} (R_{n-1}/pR_{n-1})^{[F]} \to R_{\infty}\otimes_{R_n} R_{n}/p^{1/p}R_{n} \to R_{\infty}/p^{1/p}R_{\infty} \to (R_{\infty}/pR_{\infty})^{[F]},
$$
which is obviously an $R_{\infty}$-algebra isomorphism. Hence
$$
R_{\infty} \otimes_{R_n} M_{n}^{[F]} \simeq (R_{\infty} \otimes_{R_n} M_{n+1})^{[F]}.
$$
Since the $R_{\infty}$-module $N:=R_{\infty} \otimes_{R_n} M_{n+1}$ gets mapped onto $M$ and $N$ is finitely presented, it follows that
$$
\lambda_{\infty}(M) \le \lambda_{\infty}(N)=p^d\cdot\lambda_{\infty}(N^{[F]}) \le p^d \cdot\lambda_{\infty}(N').
$$
Therefore, we get the desired formula $\lambda(M)=p^d \cdot\lambda (M^{[F]})$ for the case of finitely generated modules. 

Finally, let $ M $ be an arbitrary module in $\mathcal{M}$, and let $N \subset M$ be any finitely generated submodule. Then $N^{[F]}$ is an $R_{\infty}$-submodule of $M^{[F]}$ and $N^{[F]}$ is finitely generated over $R_{\infty}$ since the residue field of $R_{\infty}$ is assumed to be perfect. We have
$$
\lambda_{\infty}(N)=p^{d} \cdot \lambda_{\infty}(N^{[F]}) \le p^{d} \cdot \lambda_{\infty}(M^{[F]}),
$$
which gives $\lambda_{\infty}(M) \le p^{d} \cdot \lambda_{\infty}(M^{[F]})$. Conversely, let $N' \subset M^{[F]}$ be a finitely generated submodule. Then we can find a set of generators $\{v_{i}\}_{i \in I}$ of $M$, a finite subset $I'\subset I$, and an $R_{\infty}$-module $N$ such that 
$$
M^{[F]}=(\sum_{i \in I}{R_{\infty}}v_{i})^{[F]},~N=\sum_{i \in I'}{R_{\infty}}v_{i},~\mathrm{and}~N'=N^{[F]}.
$$  
Then we deduce that $\lambda_{\infty}(M) \ge \lambda_{\infty}(N)=p^{d}\cdot\lambda_{\infty}(N^{[F]})=p^{d}\cdot\lambda_{\infty}(N')$ and hence $\lambda_{\infty}(M)=p^{d}\cdot\lambda_{\infty}(M^{[F]})$, which completes the proof.
\end{proof}

In the following, we write
$$
(x^{1/{\infty}}):=\bigcup_{n \in \mathbb{N}} (x^{1/p^{n}})R_{\infty}
$$
for an element $x \in R_{\infty}$ whenever $x^{1/p^{n}} \in R_{\infty}$ for all $n \in \mathbb{N}$. Note that this is a flat $R_{\infty}$-module since it is an ascending union of principal ideals. As a corollary, one can prove:

\begin{Corollary}
Let $M \in \mathcal{M}$ such that $p^{1/p^{n}} \cdot M=0$ for all $n>0$. Then $\lambda_{\infty}(M)=0$. 
\end{Corollary}

\begin{proof}
By the definition, we easily reduce to the case in which $M \in \mathcal{M}$ is finitely generated. Then by assumption, we must have $M \simeq M^{[F]}$ as the Frobenius map on $R_{\infty}/(p^{1/{\infty}})$ is a bijection. Then we have $\lambda_{\infty}(M)=p^{d} \cdot \lambda_{\infty}(M^{[F]})=p^{d} \cdot \lambda_{\infty}(M)$ and thus $\lambda_{\infty}(M)=0$ since $\lambda_{\infty}(M)$ is finite. 
\end{proof}

Finally, we discuss the relation of $v$-almost zero modules to the vanishing of the normalized length. Let us first recall the definition of an $\fm$-adic valuation on $R$.

\begin{Definition}
Let $K$ be the field of fractions of $R$. An $\fm$-$adic$ $valuation$ on $R$ is a discrete valuation $v:K\backslash\{0\} \to \mathbb{Z}$ satisfying the property that for any non-zero $b \in R$, $v=v(b)$ is defined as the integer such that $b \in \fm^{v}$, but $b \notin \fm^{v+1}$.
\end{Definition}

If $v$ denote the $\fm$-adic valuation on $R$, then since $R_{\infty}$ is integral over $R$, it extends to a valuation on $R_{\infty}$ with value group $\mathbb{Q}$. We denote by $v$ this extended valuation for simplicity.

\begin{Proposition}
Suppose that $\lambda_{\infty}(M)=0$ for $M \in \mathcal{M}$. Then $M$ is $v$-almost zero. 
\end{Proposition}

\begin{proof}
It will suffice to consider the case in which $M$ is a cyclic $R_{\infty}$-module. Then $R_{\infty} \cdot x \simeq M$ for some $x \in M$ and there is an isomorphism $ R_{\infty}/\mathrm{Ann}_{R_{\infty}}(x) \simeq R_{\infty} \cdot x$. Now assume that there exists a non-zero $k \in \mathbb{Z}$ such that
$$
\mathrm{inf}\{ v(b)~|~ b \in \mathrm{Ann}_{R_{\infty}}(R_{\infty} \cdot x) \} \ge p^{-k}.   
$$ 
Then for any $n \in \mathbb{N}$ with $n > k$, we deduce that
$$
\ell(R_{n}\cdot x) \ge \# \{p^{\epsilon_{1}}x_{2}^{\epsilon_{2}} \cdots x_{d}^{\epsilon_{d}}~|~\sum_{i=1}^{d} \epsilon_{i} < p^{-k}~\mathrm{and}~p^{-n} \le \epsilon_{i} \}=p^{d(n-k)},
$$
where the middle term denotes the number of monomials $p^{\epsilon_{1}}x_{2}^{\epsilon_{2}} \cdots x_{d}^{\epsilon_{d}} \in R_{n} $ with specified conditions. Now it follows easily from Lemma 2.3 that
$$
\lambda_{\infty}(R_{\infty}\cdot x)=\lim_{n \to \infty}\frac{\ell(R_{n} \cdot x)}{p^{dn}} \ge \frac{1}{p^{dk}} > 0,
$$
which contradicts our assumption that $\lambda_{\infty}(M)=0$. Hence the proposition follows.
\end{proof}

\begin{Remark}
I do not know if the converse of the above proposition holds. It seems that the answer to this question is yes even if the annihilator of the module is quite complicated. Obviously, it suffices to consider the case where $M$ is a cyclic module.
\end{Remark}

\section{Applications to local cohomology modules}
\label{sec:applications-to-local-cohomology-modules}

In this section, we start with a finite extension of local domains $R \to S$ and introduce a notion of a $semi$-$perfect$ $algebra$ $S_{\infty}$ over $S$.

\begin{Definition}
Let $R \to S$ be a ring extension of domains of mixed characteristic $p>0$. Then an algebra $S_{\infty}$ over $S$ is called $semi$-$perfect$ if
\begin{enumerate}
\item
there exists a commutative square 
$$
\begin{CD}
R_{\infty} @>>> S_{\infty} \\
@AAA @AAA \\
R @>>> S 
\end{CD}
$$
such that every map is a ring extension of domains.

\item
$S_{\infty}$ is an integrally closed domain in its fraction field.

\item
the Frobenius endomorphism $F_{S_{\infty}}:S_{\infty}/pS_{\infty} \to S_{\infty}/pS_{\infty}$ is surjective.

\end{enumerate}
\end{Definition}

In what follows, $S_{\infty}$ shall denote some fixed semi-perfect algebra for $S$. Such a ring is usually obtained as a large ring extension of $S$.

\textbf{Positive characteristic}:
Let $S$ be a reduced ring of positive characteristic $p>0$. Let $S^{\infty}$ be a filtered direct limit of the system $S \hookrightarrow S^{1/p} \hookrightarrow \cdots \hookrightarrow S^{1/p^{n}} \hookrightarrow \cdots $. Then $S^{\infty}$ is the minimal perfect $S$-algebra. In fact, the Frobenius map on it is a bijection.

\textbf{Mixed characteristic}:
In this case, it is easy to see that the Frobenius map on $R_{\infty}/pR_{\infty}$ is surjective. However, since $R_{\infty}/pR_{\infty}$ is not reduced, the Frobenius map is not injective. Let $R^{+}$ be the $absolute$ $integral$ $closure$ of $R$ (see~\cite{A}); that is, it is the integral closure of $R$ in the algebraic closure of the field of fractions of $R$. Then $R^{+}$ is obviously a semi-perfect ring over any domain $S$ that is integral over $R$. In the mixed characteristic case, the notion of semi-perfect rings is quite subtle since there are not many known examples. 

\begin{Lemma}
Let $S_{\infty}$ be a semi-perfect $S$-algebra. Then there exists a ring isomorphism:
$$
\overline{F}_{S_{\infty}}:S_{\infty}/p^{1/p}S_{\infty} \simeq S_{\infty}/pS_{\infty},
$$ 
which is induced by the Frobenius endomorphism $F_{S_{\infty}}:S_{\infty}/pS_{\infty} \to S_{\infty}/pS_{\infty}$.
\end{Lemma}

\begin{proof}
Let $\overline{K}$ denote the algebraic closure of the fraction field of $S_{\infty}$. Then it is easy to see that: $F_{S_{\infty}}(\bar{x})=0$ with $x\in S_{\infty} \iff x^{p}=p \cdot \theta$ for some $\theta \in S_{\infty}$. Taking the $p$-th root of $x^{p}=p \cdot \theta$ in $\overline{K}$, we have $x=p^{1/p} \cdot \widetilde{\theta}$, where $\widetilde{\theta}$ is a root to the equation $t^{p}-\theta=0$. Hence $\widetilde{\theta} \in S_{\infty}[1/p]$, which gives $\widetilde{\theta} \in S_{\infty}$, since $S_{\infty}$ is assumed integrally closed. This proves the desired claim.
\end{proof}

In what follows, we let $v$ be a valuation on $R_{\infty}$ that extends an $\fm$-adic valuation of $R$, and let $d=\mathrm{dim}~R$. A ring isomorphism $\overline{F}_{S_{\infty}}:S_{\infty}/p^{1/p}S_{\infty} \simeq S_{\infty}/pS_{\infty}$ yields a well-defined homomorphism on $\mathrm{\check{C}}$ech complexes:
$$
\overline{F}_{S_{\infty}*}:C^{\bullet}(S_{\infty}/p^{1/p}S_{\infty}) \simeq C^{\bullet}(S_{\infty}/pS_{\infty})^{[F]}
$$
and the induced map on local cohomology modules:
$$
\overline{F}_{S_{\infty}*}:H^{i}_{\fm}(S_{\infty}/p^{1/p}S_{\infty}) \simeq H^{i}_{\fm}(S_{\infty}/pS_{\infty})^{[F]}
$$
for $i \in \mathbb{N}$. Now let us suppose the following condition: For $\epsilon \in \mathbb{Q}$ with $\epsilon > 0$, there exists an element $r \in R_{\infty}$ such that
$$
\lambda_{\infty}(r \cdot H^{i}_{\fm}(S_{\infty})) < \infty \cdots(*)
$$
with $v(r)< \epsilon$ for some $i < d$. Recall that the condition $(*)$ simply says that the local cohomology $H^{i}_{\fm}(S_{\infty})$ has almost finite length with respect to the valuation $v$. In order to prove the main theorem, we need the following lemma taken from~\cite{GR2}:

\begin{Lemma}
Let $0 \to N' \to N \to N'' \to 0$ be a short exact sequence in the category $\mathcal{M}$ and let $a, b \in R_{\infty}$. Then the following hold:
\begin{enumerate}

\item[$\mathrm{(1)}$]
$\lambda_{\infty}(ab N) \le \lambda_{\infty}(a N') + \lambda_{\infty}(b N'')$.

\item[$\mathrm{(2)}$]
Suppose that $\lambda_{\infty}(N')=0$ $($resp.~$\lambda_{\infty}(N'')=0$$)$. Then $\lambda_{\infty}(a N)=\lambda_{\infty}(a N'')$ $($resp.~$\lambda_{\infty}(a N)=\lambda_{\infty}(a N')$$)$.
\end{enumerate}
\end{Lemma}

\begin{proof}
(1): First of all, we notice that there is a commutative diagram:
$$
\begin{CD}
0 @>>> bN \cap N' @>>> bN @>>> bN'' @>>> 0 \\
@. @ VaVV @ VaVV @ VaVV \\
0 @>>> bN \cap N' @>>> bN @>>> bN'' @>>> 0 \\
\end{CD}
$$
where the horizontal sequences are short exact. It is easy to show that there is a surjection $bN''/L \twoheadrightarrow abN''$ where $L$ denotes the image of the kernel of the multiplication map $a: bN \to bN$ under the surjection $bN \twoheadrightarrow bN''$. Hence there follows an exact sequence
$$
0 \to a(bN \cap N') \to  abN \to bN''/L \to 0.
$$
So $\lambda_{\infty}(abN)=\lambda_{\infty}(a(bN \cap N'))+\lambda_{\infty}(bN''/L) \le \lambda_{\infty}(aN')+\lambda_{\infty}(bN'')$. Hence part (1) is proved.

(2): This follows easily from the short exact sequence
$$
0 \to aN \cap N' \to aN \to aN'' \to 0~(\mathrm{resp}.~0 \to aN' \to aN \to P \to 0)
$$
where $P$ is a quotient of $N''$, and the fact that $\lambda_{\infty}(aN \cap N')=0$ (resp. $\lambda_{\infty}(P)=0$).
\end{proof}

Now we are ready to prove the almost vanishing theorem for local cohomology modules. Let $M$ be a module over a ring $A$ and let $a \in A$. Then we use the following notation: $M_{a}:=\{m \in M~|~a m=0\}$.

\begin{Theorem}
Let $R \to S$ be a module-finite extension of domains and we fix $k \in \mathbb{N}$ with $0<k<d$. Assume the following conditions:
\begin{enumerate}

\item[$\mathrm{(1)}$]
$S_{\infty}$ is a semi-perfect algebra over $S$.

\item[$\mathrm{(2)}$]
$\lambda_{\infty}(H^{k-1}_{\fm}(S_{\infty}))=0$.

\item[$\mathrm{(3)}$]
$H^{k}_{\fm}(S_{\infty})$ has almost finite length.
\end{enumerate}
Then the local cohomology module $H^{k}_{\fm}(S_{\infty})$ is $v$-almost zero.
\end{Theorem}

\begin{proof}
Notice that by assumption we have a well-defined $R_{\infty}$-isomorphism:
$$
\overline{F}_{S_{\infty}*}:H^{i}_{\fm}(S_{\infty}/p^{1/p}S_{\infty}) \simeq H^{i}_{\fm}(S_{\infty}/pS_{\infty})^{[F]}.
$$
Note that since $S_{\infty}$ is a domain, $H^{0}_{\fm}(S_{\infty})=0$. The following diagram of standard short exact sequences:
$$ 
\begin{CD}
0 @>>> S_{\infty} @>p^{1/p}>> S_{\infty} @>>> S_{\infty}/p^{1/p}S_{\infty} @>>> 0 \\
@. @. @. @V \wr V\overline{F}_{S_{\infty}}V \\
0 @>>> S_{\infty} @>p>> S_{\infty} @>>> S_{\infty}/pS_{\infty} @>>> 0 
\end{CD}
$$
yields a diagram of local cohomology modules:
$$
\begin{CD}
\cdots @>>> H^{k-1}_{\fm}(S_{\infty}/p^{1/p}S_{\infty}) @>>> H^{k}_{\fm}(S_{\infty}) @>p^{1/p}>> H^{k}_{\fm}(S_{\infty}) @>>> \cdots \\
@. @ V \wr V \overline{F}_{S_{\infty}*} V \\
\cdots @>>> H^{k-1}_{\fm}(S_{\infty}/pS_{\infty}) @>>> H^{k}_{\fm}(S_{\infty}) @>p>> H^{k}_{\fm}(S_{\infty}) @>>> \cdots \\
\end{CD}
$$
Let us fix a positive $\epsilon \in \mathbb{Q}$. Then there exists $r \in R_{\infty}$ with $v(r) < \epsilon$ such that $\lambda_{\infty}(r \cdot H^{k}_{\fm}(S_{\infty})) < \infty$. From the above exact sequence there follows a surjection: 
$$
H^{k-1}_{\fm}(S_{\infty}/p^{1/p}S_{\infty}) \twoheadrightarrow H^{k}_{\fm}(S_{\infty})_{p^{1/p}}~(\mathrm{resp}.~H^{k-1}_{\fm}(S_{\infty}/pS_{\infty}) \twoheadrightarrow H^{k}_{\fm}(S_{\infty})_{p})
$$
Then by the Frobenius pull-back formula, we deduce
$$
\lambda_{\infty}(r \cdot H^{k-1}_{\fm}(S_{\infty}/p^{1/p}S_{\infty}))=
\lambda_{\infty}(r \cdot (H^{k-1}_{\fm}(S_{\infty}/pS_{\infty}))^{[F]})=\frac{1}{p^d} \cdot \lambda_{\infty}(r^{p} \cdot H^{k-1}_{\fm}(S_{\infty}/pS_{\infty}))<\infty,
$$
which, together with Lemma 3.3 (2), yields the equality:
$$
\lambda_{\infty}(r \cdot H^{k}_{\fm}(S_{\infty})_{p^{1/p}})=\frac{1}{p^d} \cdot\lambda_{\infty}(r^{p} \cdot H^{k}_{\fm}(S_{\infty})_{p})\cdots(1)
$$
To ease the notation, we set $N:=H^{k}_{\fm}(S_{\infty})$. Let us consider the filtration of $N_{p}$ as follows:
$$
N_{p^{1/p}} \subset N_{p^{2/p}} \subset \cdots \subset N_{p^{t/p}} \subset \cdots \subset N_{p}.
$$
Then since $p^{(t-1)/p} N_{p^{t/p}} \subset N_{p^{1/p}}$, the multiplication map $p^{(t-1)/p}:N_{p^{t/p}} \to N_{p^{1/p}}$ induces an injective map $N_{p^{t/p}}/N_{p^{(t-1)/p}} \hookrightarrow N_{p^{1/p}}$ and we have
$$
\lambda_{\infty}(N_{p^{t/p}}/N_{p^{(t-1)/p}}) \le \lambda_{\infty}(N_{p^{1/p}})\cdots(2).
$$
Let us next consider the short exact sequence
$$
0 \to N_{p^{(t-1)/p}} \to N_{p^{t/p}} \to N_{p^{t/p}}/N_{p^{(t-1)/p}} \to 0.
$$
Then from Lemma 3.3 we have 
$$
\lambda_{\infty}(r^{t} (N_{p^{t/p}})) \le \lambda_{\infty}(r^{t-1} (N_{p^{(t-1)/p}}))+\lambda_{\infty}(r (N_{p^{t/p}}/N_{p^{(t-1)/p}}))
$$
for $t=2,\ldots,p$. An induction on $t$ yields 
$$
\lambda_{\infty}(r^{p} (N_{p})) \le \sum_{t=1}^{p} \lambda_{\infty}(r (N_{p^{t/p}}/N_{p^{(t-1)/p}})) \cdots(3),
$$
and it follows easily from $(2)$ and $(3)$ that
$$
\lambda_{\infty}(r^{p} (N_{p})) \le p \cdot \lambda_{\infty}(r (N_{p^{1/p}})).
$$
We can now deduce from this inequality and $(1)$:
$$
\frac{\lambda_{\infty}(r^{p} (N_{p}))}{p^{d}}=\lambda_{\infty}(r (N_{p^{1/p}})) \ge \frac{\lambda_{\infty}(r^{p} (N_{p}))}{p},
$$
which is possible only when $\lambda_{\infty}(r (N_{p^{1/p}}))=0$ (or equivalently, $\lambda_{\infty}(r^{p} (N_{p}))=0$). For simplicity, we denote $r^{p}$ by $r$, so that $\lambda_{\infty}(r(N_{p}))=0$. Now let $\eta \in N$. Since $N$ is $p$-torsion, we have $\eta \in N_{p^{n}}$ for sufficiently large $n \in \mathbb{N}$. Then we may use the fact that $p \cdot N_{p^t} \subset N_{p^{t-1}}$, Lemma 3.3 together with the following short exact sequence inductively
$$
\begin{CD}
0 @>>> N_{p} @>>> N_{p^t} @>p>> p \cdot N_{p^t} @>>> 0
\end{CD}
$$
to deduce that
$$
\lambda_{\infty}(r^{n} (N_{p^n}))=0.
$$ 
In light of Proposition 2.15, for $\epsilon \in \mathbb{Q}$ as above, there exists an element $s \in R_{\infty}$ such that $(s r^{n}) \eta=0$ with $v(s r^{n})=v(s)+v(r^{n})<\epsilon+n\epsilon=\epsilon(n+1)$. If we have chosen $\epsilon$ sufficiently small, we can make $\epsilon(n+1)$ also sufficiently small, as $n \in \mathbb{N}$ depends on the choice of $\eta \in N$. Hence this proves the desired claim.
\end{proof}

\begin{Remark}
An $R_{\infty}$-module $M$ with $\lambda_{\infty}(M)=0$ is $v$-almost zero. However, it might be the case that the annihilator of every element of $M$ is very complicated and $M$ is not finitely generated. It looks reasonable to ask whether the local cohomology module of any semi-perfect ring is $v$-almost zero in light of the above theorem. For some relevant results of the annihilator, see also \cite{R1},~\cite{Sc}. 
\end{Remark}

Now we remind the reader that $R_{n}$ and $d=\mathrm{dim}~R$ are defined as previously. Let $R_{n} \to S_{n}$ be a module-finite extension of domains with $p^{1/p} \in S_{n}$. Let $S_{\infty}$ be any semi-perfect algebra over $S_{n}$ and let $\phi:H^{k}_{\fm}(S_{n}) \to H^{k}_{\fm}(S_{\infty})$ be the map induced by the inclusion $S_{n} \to S_{\infty}$. Then define the map 
$$
\Phi:H^{k}_{\fm}(S_{n}) \otimes_{R_n} R_{\infty} \to H^{k}_{\fm}(S_{\infty})
$$
by $\Phi(a \otimes b):=\phi(a)b$ for $a \in H^{k}_{\fm}(S_{n})$ and $b \in R_{\infty}$. We also denote by  
$$
\Phi_{p^{1/p}}:H^{k}_{\fm}(S_{n})_{p^{1/p}} \otimes_{R_n} R_{\infty} \to H^{k}_{\fm}(S_{\infty})_{p^{1/p}}
$$
the restriction map of $\Phi$ to the submodule annihilated by $p^{1/p}$. Alternatively, one may regard the map $\Phi$ to be induced by the natural map $S_{n} \otimes_{R_n} R_{\infty} \to S_{\infty}$ as $R_{\infty}$ is flat over $R_{n}$.

\begin{Theorem}
Keeping the notations and hypothesis as above, assume that we have $H^{k-1}_{\fm}(S_{n})=0$ and $\ell(H^{k}_{\fm}(S_{n}))<\infty$ for $0<k<d$. Then
$$
\lambda_{\infty}(\mathrm{Im}(\Phi_{p^{1/p}})) \le \frac{1}{p^{d(n+1)-1}} \cdot \ell(H^{k}_{\fm}(S_{n})_{p^{1/p}}).
$$
In particular, $\Phi$ is not injective if $H^{k}_{\fm}(S_{n})$ is non-zero.
\end{Theorem}

\begin{proof}
First, we note that $H^{k}_{\fm}(S_{n}) \ne 0 \iff H^{k}_{\fm}(S_{n})_{p^{1/p}} \ne 0$. Since $S_{n}$ is a domain, one has the following short exact sequence:
$$
\begin{CD}
0 @>>> S_{n} @>b>> S_{n} @>>> S_{n}/bS_{n} @>>> 0 \\
\end{CD}
$$
for any non-zero $b \in S_{n}$. Then by assumption, the associated long exact sequence
$$
\begin{CD}
0 @>>> H^{k-1}_{\fm}(S_{n}/bS_{n}) @>>> H^{k}_{\fm}(S_{n}) @>b>> H^{k}_{\fm}(S_{n}) @>>> \cdots \\
\end{CD}
$$
yields an isomorphism $H^{k}_{\fm}(S_{n})_{b} \simeq H^{k-1}_{\fm}(S_{n}/bS_{n}) \cdots (1)$. Let us consider the following commutative square of $R_{\infty}$-modules:
$$
\begin{CD}
S_{n}/p^{1/p}S_{n} \otimes_{R_n} R_{\infty} @>>> S_{\infty}/p^{1/p}S_{\infty}\\
@V \overline{F}_{S_n} VV @V \overline{F}_{S_{\infty}} V\wr V \\
(S_{n}/pS_{n}\otimes_{R_n} R_{\infty})^{[F]} @>>> (S_{\infty}/pS_{\infty})^{[F]} \\\end{CD}
$$
in which both of the horizontal maps are induced by the natural map $S_{n}/p^{1/p}S_{n} \to S_{\infty}/p^{1/p}S_{\infty}$ (resp. $S_{n}/pS_{n} \to S_{\infty}/pS_{\infty}$) and the vertical maps are induced by the Frobenius map. Then by use of the isomorphism $(1)$ together with the flatness of $R_{\infty}$ over $R_{n}$, we have the induced map on local cohomology modules:
$$
\begin{CD}
H^{k}_{\fm}(S_{n})_{p^{1/p}}\otimes_{R_n} R_{\infty} @ >\widetilde{\Phi}_{p^{1/p}} >> H^{k-1}_{\fm}(S_{\infty}/p^{1/p}S_{\infty}) \\
@V \overline{F}_{S_{n}*} VV @V \overline{F}_{S_{\infty}*} V \wr V \cdots(2)\\
(H^{k}_{\fm}(S_{n})_{p}\otimes_{R_n} R_{\infty})^{[F]} @ > \widetilde{\Phi}_{p} >> H^{k-1}_{\fm}(S_{\infty}/pS_{\infty})^{[F]} \\
\end{CD}
$$
By assumption, the module $H^{k}_{\fm}(S_{n})_{p^{1/p}}\otimes_{R_n} R_{\infty}$ is finitely presented in $\mathcal{M}$, and we get
$$
\lambda_{\infty}(H^{k}_{\fm}(S_{n})_{p^{1/p}}\otimes_{R_n} R_{\infty})=\frac{1}{p^{dn}}\cdot \ell(H^{k}_{\fm}(S_{n})_{p^{1/p}}).
$$
and similarly
$$
\lambda_{\infty}(H^{k}_{\fm}(S_{n})_{p}\otimes_{R_n} R_{\infty})=\frac{1}{p^{dn}}\cdot \ell(H^{k}_{\fm}(S_{n})_{p}).
$$
On the other hand, an easy calculation shows that $\ell(H^{k}_{\fm}(S_{n})_{p}) \le p \cdot \ell (H^{k}_{\fm}(S_{n})_{p^{1/p}})$. Hence it follows from the diagram $(2)$ and the Frobenius pull-back formula that
$$
\lambda_{\infty}(\mathrm{Im}(\widetilde{\Phi}_{p^{1/p}})) \le \lambda_{\infty}((\mathrm{Im}(\widetilde{\Phi}_{p}))^{[F]})=\frac{1}{p^d} \cdot \lambda_{\infty}(\mathrm{Im}(\widetilde{\Phi}_{p})) \le \frac{1}{p^{d(n+1)}} \cdot \ell(H^{k}_{\fm}(S_{n})_{p})
$$
$$
\le \frac{p}{p^{d(n+1)}} \cdot \ell(H^{k}_{\fm}(S_{n})_{p^{1/p}})=\frac{1}{p^{d(n+1)-1}} \cdot \ell(H^{k}_{\fm}(S_{n})_{p^{1/p}}).
$$
Let $j_{b}:H^{k-1}_{\fm}(S_{\infty}/bS_{\infty}) \twoheadrightarrow H^{k}_{\fm}(S_{\infty})_{b}$ denote the surjection map induced by the cohomology exact sequence associated to:
$$
\begin{CD}
0 @>>> S_{\infty} @>b>> S_{\infty} @>>> S_{\infty}/bS_{\infty} @>>> 0
\end{CD}
$$
for any non-zero $b \in S_{\infty}$. Then by assumption, it is easy to see that $\Phi_{b}=j_{b} \circ \widetilde{\Phi}_{b}$, from which $\lambda_{\infty}(\mathrm{Im}(\widetilde{\Phi}_{p^{1/p}})) \ge \lambda_{\infty}(\mathrm{Im}(\Phi_{p^{1/p}}))$. Hence the theorem follows easily from this inequality.

The non-injectivity of the map $\Phi$ is due to the following equality:
$$
\lambda_{\infty}(H^{k}_{\fm}(S_{n})_{p^{1/p}} \otimes_{R_{n}} R_{\infty})=\frac{1}{p^{dn}} \cdot \ell(H^{k}_{\fm}(S_{n})_{p^{1/p}}).
$$
\end{proof}

The next corollary is a weak version of Heitmann's Direct Summand Theorem.

\begin{Corollary}
Let $R_{n} \to S_{n}$ be a module-finite extension of normal domains and $p^{1/p} \in S_{n}$, and let $d>2$. Then the map
$$
\Phi:H^{2}_{\fm}(S_{n}) \otimes_{R_n} R_{\infty} \to H^{2}_{\fm}(S^{+})
$$
is not injective if $H^{2}_{\fm}(S_{n})$ is non-zero.
\end{Corollary}

\begin{proof}
By Serre's normality criterion, it follows that $H^{1}_{\fm}(S_{n})=0$ and $\ell(H^{2}_{\fm}(S_{n}))<\infty$ by local duality. Then the corollary follows from Theorem 3.6.
\end{proof}

There is another implication for splinters studied by A. Singh. Recall that a Noetherian domain $A$ is a $splinter$ if $A$ is a direct summand, as an $A$-module, of every module-finite extension domain. It is easy to see that a splinter is a normal domain.

\begin{Corollary}
In addition to the hypothesis of Theorem 3.6, assume that $S \otimes_{R} R_{\infty}$ is a domain and $H^{k}_{\fm}(S)$ is non-zero. Then there exists a flat module-finite extension $S \to T$ such that $T$ is not a splinter. \end{Corollary}

\begin{proof}
Since $\Phi:H^{k}_{\fm}(S) \otimes_{R} R_{\infty} \to H^{k}_{\fm}(S^{+})$ is not injective and $S \otimes_{R} R_{\infty}=\bigcup_{n \ge 0} S \otimes_{R} R_{n}$, the map $H^{k}_{\fm}(S) \otimes_{R} R_{n} \to H^{k}_{\fm}(S^{+})$ is not injective for some $n \ge 0$. Now it is easy to see that the domain $T:=S \otimes_{R} R_{n}$ is flat over $S$ and is not a splinter.
\end{proof}

Finally, we end this section with an example of a normal domain that is not a splinter, due to A. Singh.

\begin{Example}
Let 
$$
S:=\frac{\mathbb{Z}_{2}[[x,y]]}{(x^{k}-y^{l}-2^{m})}
$$
where $\mathbb{Z}_{2}$ is the ring of 2-adic integers and $(k,l,m) > (2,2,2)$. Then $S$ is shown to be a normal domain. Now take $u:=\sqrt{x^{k-2}y^{l-2}} \notin S$, and $v:=2^{-1}(x^{k-1}-uy) \notin S$. Then $u$ is integral over $S$ and it is easy to verify that $v$ is subject to the following equation:
$$
v^{2}+uyv-x^{k-2}2^{m-2}=0
$$
Hence $S \to T:=S[u,v]$ is module-finite. Then since $x^{k-1}=uy+v2$, if $S$ were a direct summand of $T$, we would have that $x \in (y,2)T$, while $x \notin (y,2)S=(y,2)T \cap S$, which is a contradiction.
\end{Example}

\textbf{Acknowledgments} \\
I would like to express my sincere gratitude to P. Roberts for explaining to me lots of new ideas to start this research project. I also want to extend my gratitude to G. Piepmeyer for allowing me to look at his notes, which led to many corrections and improvements. My final gratitude goes to A. Singh and the referee, who kindly provided me with useful comments.

\end{document}